\documentclass[12pt,leqno, a4paper]{article}
\usepackage[latin1]{inputenc}
\usepackage[english]{babel}
\renewcommand{\theequation}{\thesection.\arabic{equation}}
\newtheorem{theorem}{Theorem}[section]
\newtheorem{lemma}[theorem]{Lemma}
\newtheorem{prop}[theorem]{Proposition}

\newtheorem{fact}[theorem]{Fact}

\newtheorem{remark}[theorem]{Remark}

\newcommand{\eqnsection}{
\renewcommand{\theequation}{\thesection.\arabic{equation}}
        \makeatletter
        \csname  @addtoreset\endcsname{equation}{section}
        \makeatother}
\eqnsection

\newcommand{\bp}{{\bf P}}
\newcommand{\br}{I\!\!R}
\newcommand{\be}{{\bf E}}

\def \bt{\begin{theorem}}
\def \et{\end{theorem}}
\def \bea{\begin{eqnarray}}
\def \eea{\end{eqnarray}}
\def \bas{\begin{eqnarray*}}
\def \eas{\end{eqnarray*}}

\def \al{\alpha}
\def \bb{\beta}
\def \ga{\gamma}

\def \la{\lambda}

\def \si{\sigma}

\def \ze{\zeta}

\def \ff{\infty}

\def \wt{\widetilde}

\def \({\left(}
\def \){\right)}

\def \nn{\nonumber}

\def \bc{\begin{center} }
\def \ec{\end{center} }
\def \bs{\begin{slide} }
\def \es{\end{slide} }

\def\square{{\vcenter{\vbox{\hrule height.3pt
                \hbox{\vrule width.3pt height5pt \kern5pt
                   \vrule width.3pt}
                \hrule height.3pt}}}}

\begin{document}


\title{{\bf \LARGE  Large time asymptotics for the density of a branching
Wiener process}}

\author{P\'{a}l R\'{e}v\'{e}sz, \, Jay Rosen\thanks {Research supported,
in part, by grants from the NSF  and from PSC-CUNY.} \,\, and
\, Zhan Shi\\ {\small\it Technische Universit\"at Wien, City University of
New York $\; \&$ Universit\'e Paris VI}}
\date{}

\maketitle

\begin{abstract} Given an $\br^d$-valued supercritical branching
Wiener process, let $\psi(A,T)$ be the number of particles in
$A\subset\br^d$ at time $T,\,\ (T=0,1,2,\ldots)$. We provide a complete
asymptotic expansion of $\psi(A,T)$ as
$T\rightarrow\infty$, generalizing the work of X.~Chen (\cite{chen}).
\end{abstract}




\section{Introduction}
\label{s:intro}

Consider the following model in $\br^d$ (with $d\ge 1$):
\begin{itemize}
\item[(i)] a particle starts from the origin in $\br^d$ and executes a
               Wiener process $W(t)\in\br^d$,
\item[(ii)] arriving at time $t=1$ at the new location $W(1)$, it dies,
\item[(iii)] at death it is replaced by $Y$ of\/fspring where
\begin{eqnarray*}
     &&\bp\{Y=\ell\}=p_\ell,\quad (\ell=0,1,2,\ldots)
        \\
     &&1<\sum_{\ell=0}^\infty \ell p_\ell=m<\infty,
        \\
     &&0<\sum_{\ell=0}^\infty(\ell-m)^2p_\ell=\sigma^2<\infty,
\end{eqnarray*}
\item[(iv)] each of\/fspring, starting from where its ancestor dies,
                executes a Wiener process (from its starting point) and repeats
                the above given steps and so on. All Wiener processes and
                of\/fspring-numbers are assumed independent of each other.
\end{itemize}

Let
$$
\lambda(x,t)=\left\{\begin{array}{ll} 1&{\rm if\ }x\in\br^d\ {\rm is\
occupied\ by\ a\ particle\ at\ time\ } t,\\ 0 & {\rm otherwise}.
\end{array}\right.
$$

\noindent We write
$$
\psi(A,t)=\sum_{x\in A}\lambda(x,t),
$$

\noindent which stands for the number of particles at time $t$ located at
$A\subset\br^d$. In particular, $\psi(\br^d,t)$ is the total number of
particles alive at time $t$.

Since the branching is supercritical, it is well-known (Athreya and Ney
\cite{athreya-ney}, p.~9) that
\begin{equation}
       N_0
     :=\lim_{T\rightarrow\infty}{\frac{\psi(\br^d,T)}{m^T}}\quad {\rm
a.s.},
       \label{N0}
\end{equation}

\noindent exists (and is finite), and that $\bp(N_0>0)>0$.

The limit properties of $\psi(A,T),\,T\rightarrow\infty$, were studied by
Chen (\cite{chen}) who proved

\medskip

\renewcommand{\thetheorem}{\Alph{theorem}.}
\setcounter{theorem}{0}
\begin{theorem}
        There exist random variables $N_1$ and $N_2$ ($N_1$ being
        $\br^d$-valued) such that for any Borel set
        $A\subset \br^d$ with $\int_A \|x\|^2 dx <\infty$, we have, almost
        surely when $T\to \infty$,
$$  (2\pi T)^{d/2}{\frac{\psi(A,T)}{m^T}} =N_0\int_A dx -{\frac{1}{2T}}
\int_A(N_0 \|x\|^2 -2N_1\cdot x + N_2)dx +o(T^{-1}).
$$
\end{theorem}

\medskip

This result plays an important role in R\'ev\'esz (\cite{revesz04}) in the
study of the concentration of particles in the branching process.

The goal of this paper is to provide a complete asymptotic expansion for
$\psi(A,T)/m^T$ as $T\rightarrow\infty$. Let us first introduce some
notation.

If $\alpha=( \al_{ 1},\ldots,\al_{ d})\in Z_{ +}^{ d }$ and
$x=(x_{ 1},\ldots,x_{ d})\in R^{ d }$ we use the notation
$|\alpha|= \al_{ 1}+\cdots+\al_{ d}$, $\,\al !=\prod_{ i=1}^{ d} \al_{ i}!$,
$\,x^\al=\prod_{ i=1}^{ d}x^{ \al_{ i}}_{ i}$ and
\begin{equation} M_{ \al}(A)=\int_A x^\al\,dx. \label{m.1}
\end{equation} If also $\bb\in Z_{ +}^{ d }$ we will write
$\bb\preceq\al$ to mean that
$\bb_{ i}\leq\al_{ i}$ for all $i$, and if $\bb\preceq\al$ we set
\begin{equation} { \al\choose\bb }=\prod_{ i=1}^{ d}{ \al_{
i}\choose\bb_{ i} }.
\label{m.2}
\end{equation}

\noindent Here is the main result of the paper:

\medskip

\renewcommand{\thetheorem}{\thesection.\arabic{theorem}}
\setcounter{theorem}{0}
\begin{theorem}
\label{t:main}
        There exist random variables $(N_{ \al}, \; \alpha\in Z_{ +}^{ d })$
        such that for any $k\ge 1$ and any bounded Borel set
$A\subset\br^d$,
        when $T\to \infty$,
\bea &&\qquad
        (2\pi T)^{d/2}{\frac{\psi(A,T)}{m^T}}    \label{Q}\\ &&
=\sum_{n=0}^k {(-T)^{-n} \over 2^{n} }\sum_{|\alpha|=n}{ 1\over \al !}
\sum_{\bb\preceq 2\al} { 2\al\choose\bb }   ( -1)^{ |\bb|}
M_{\bb}(A)N_{ 2\al-\bb} +o(T^{-k}),
        \qquad \hbox{\rm a.s.}\nn
\eea
\end{theorem}

\begin{remark} {\rm The random variables $(N_{ \al}, \; \alpha\in Z_{
+}^{ d })$ are described in the proof of Theorem \ref{t:main}. They are
limits of explicit martingales related to the branching Wiener process. }
\end{remark}

Although the distributions of the random variables $(N_{ \al}, \;
\alpha\in Z_{ +}^{ d })$ are not known, Theorem  \ref{t:main} can
nevertheless be used to make predictions to any degree of accuracy.

To see this, choose an integer $k$ and disjoint sets $(A_{ \al}\subseteq
\br^d,
\; |\alpha|\leq k)$.
Consider (1.4) for each $A_\alpha$. Then we have a linear system of equations with the unknowns $N_{2\alpha-\beta}$. One can solve this system
of equations if the corresponding determinant is not equal to $0$. It is easy to see that we can choose the sets $A_\alpha$ such that the
determinant is not $0$ for any $T\ (T=1,2,\ldots)$. 
Observe the number of particles of a branching
Wiener process which are located in the above given sets $(A_{ \al}, \; |\alpha|\leq k)$ at time
$T_0$. Having  these observations one can evaluate the actual values of
the random variables $(N_{ \al}, \; |\alpha|\leq k)$ with an error term
$o(T_0^{-k})$.  Having these values one can use Theorem
\ref{t:main} to get the values of the process $(2\pi
T)^{d/2}\psi(A,T)/m^T$ for any
   $A\subseteq \br^d,\ T\geq T_0$ with an error term $o(T_0^{-k})$.

\medskip

The proof of Theorem \ref{t:main} is presented in Section \ref{s:proof}.
In Section \ref{s:Lp} we show that if the offspring distribution $Y$ has
$p$  moments for some even integer $p$ then the martingales described
in Remark 1.2 converge to the random variables $(N_{ \al}, \;
\alpha\in Z_{ +}^{ d })$ in
$L^{ p}$.

\section{The proof}
\label{s:proof}

We start with a preliminary result concerning the transition kernel of the
Wiener process. Let
$$ p^{ ( d)}_t(x) = {1\over (2\pi t)^{d/2}} \exp\left( -  {\|x\|^2\over
2t}\right).
$$

Define the Hermite polynomials by
\begin{equation} H_{ n}( x,t)=\sum_{ j=0}^{ [n/2]}{ n!\over j! (
n-2j)!}\({-t
\over 2}\)^{ j}x^{ n-2j}.\label{1.1}
\end{equation}

\medskip

\begin{lemma}
\label{l:kernel}
        For any $0<t<T$ and any $x\in \br^1$,
\begin{equation}
        p^{ ( 1)}_{T-t}(x)=  {1\over (2\pi T)^{1/2}} \sum_{n=0}^\infty
{(-T)^{-n}
        \over 2^{n} n!} H_{2n} (x, t).
        \label{kernel1}
\end{equation}
\end{lemma}

\medskip

\noindent {\it Proof.} Let us recall the Hermite polynomials:
\begin{eqnarray}
        H_n(x)
     &=&  (-1)^n e^{x^2} {d^n\over d x^n} (e^{x^2})
        \nonumber
        \\
     &=& n! \sum_{j=0}^{\lfloor n/2\rfloor} {(-1)^j\over j! (n-2j)!}
        (2x)^{n-2j}, \qquad x\in \br,
        \label{Hermite}
\end{eqnarray} so that
\begin{equation}
      H_n(x, t) =(t/2)^{n/2} H_n\left( {x\over \sqrt{2t}}\right) , \qquad
x\in
\br,
       \; \; t>0.
       \label{Hermite(.,.)}
\end{equation}

We use the following identity, see for example Lebedev (\cite{lebedev},
p.~75): for any $a>0$ and $y\in \br$,
$$ e^{-a^2y^2} = \sum_{n=0}^\infty {(-1)^n a^{2n} \over 2^{2n} n!
(1+a^2)^{n+(1/2)}} H_{2n} (y).
$$

\noindent Taking $y=x/\sqrt{2t}\, \in \br^1$ and $a= \sqrt{t/(T-t)}$,
and multiplying both sides by $(2\pi (T-t))^{-1/2}$, we readily get
(\ref{kernel1}).
\hfill$\diamondsuit$

\bigskip

If $\alpha=( \al_{ 1},\ldots,\al_{ d})\in Z_{ +}^{ d }$ and
$x=(x_{ 1},\ldots,x_{ d})\in R^{ d }$ we use the notation
\begin{equation}
       H_{\al} (x, t) =\prod_{ i=1}^{ d}H_{\al_{i}} (x_{ i}, t).
       \label{2k.1}
\end{equation}

\medskip

\begin{lemma}
\label{ld:kernel}
        For any $0<t<T$ and any $(x,y)\in \br^d \times \br^d$,
\begin{equation}
        p^{ ( d)}_{T-t}(x)=  {1\over (2\pi T)^{d/2}} \sum_{n=0}^\infty
{(-T)^{-n}
        \over 2^{n} }\sum_{|\alpha|=n}{ 1\over \al !}
        H_{2\al} (x, t) ,
        \label{kernel1d}
\end{equation} and
\bea &&
        p^{ ( d)}_{T-t}(x-y)    \label{mkernel1d} \\ &&=  {1\over (2\pi
T)^{d/2}}
\sum_{n=0}^\infty {(-T)^{-n}
        \over 2^{n} }\sum_{|\alpha|=n}{ 1\over \al !}
\sum_{\bb\preceq 2\al} { 2\al\choose\bb } (-x)^{\bb} H_{2\al-\bb} (y,
t).
\nn
\eea
\end{lemma}

\medskip

\noindent {\it Proof.} Since for $x=(x_{1},\ldots,x_{ d} )\in \br^d$
\begin{equation} p^{ ( d)}_t(x)=\prod_{ i=1}^{ d}p^{ ( 1)}_t(x_{
i}),\label{m2.1}
\end{equation} (\ref{kernel1d}) follows from (\ref{kernel1}). To obtain
(\ref{mkernel1d}) we use the fact that
\begin{equation} H_{n} (x+y, t)=\sum_{ j=0}^{ n}{n\choose j } x^{ n-j}
H_{j}(y, t).\label{m2.2}
\end{equation}

For this we recall that (Lebedev \cite{lebedev}, p.~60)
\begin{equation}
\sum_{ n=0}^{ \ff}{s^{ n} \over n!}H_{ n}( x)=e^{ 2s x-s^{ 2}}\label{3m.1}
\end{equation} so that
\begin{equation}
\sum_{ n=0}^{ \ff}{s^{ n} \over n!}H_{ n}( x,t)=e^{ s x-ts^{
2}/2}.\label{3m.2}
\end{equation} Then
\bea &&
\sum_{ n=0}^{ \ff}{s^{ n} \over n!}H_{ n}( x+y,t) =e^{ s (x+y)-ts^{2}/2}
\label{3m.222}
\\ &&=e^{ s x}e^{ s y-ts^{ 2}/2}=\sum_{ k=0}^{ \ff}{s^{ k}x^{ k}
\over k!}
\sum_{ j=0}^{ \ff}{s^{ j} \over j!}H_{ j}( y,t) , \nn
\eea and comparing powers of  $s^{ n}$ proves
(\ref{m2.2}).\hfill$\diamondsuit$

\bigskip

Now we turn to the study of the branching Wiener process. Clearly, for
any
$T\ge 1$ and $A\subset \br^d$,
$$
\be\left( \psi(A,T) \, | \, {\cal F}(T-1)\right) = m\int_A\sum_y p^{  (
d)}_1(y-x) \lambda(y,T-1) dx ,
$$

\noindent (as usual, ${\cal F}(t)$ denoting the $\sigma$-algebra induced
by the branching process until time $t$). A simple argument by
induction yields that for all $0<t<T$,
\begin{equation}
        \be\left( \psi(A,T) \, | \, {\cal F}(t)\right) = m^{T-t}\int_A\sum_y
        p^{ ( d)}_{T-t}(y-x) \lambda(y,t) dx.
        \label{E(psi|F)}
\end{equation}

\noindent It turns out that $\psi(A,T)$ is quite close to its conditional
expectation, as is confirmed by the following results.
\medskip

\begin{fact} {\rm (R\'ev\'esz \cite{revesz}, (6.16))}
\label{f:cond-exp}
        Fix $\gamma\in (0,1)$ and let $t= \lfloor T^\gamma\rfloor$. Let
        $A\subset \br^d$ be a bounded Borel set. Let $\varepsilon>0$. We
        have, almost surely for $T\to \infty$,
\begin{equation}
        { \frac{\psi(A,T)}{m^T}} - {\frac{1}{m^t}} \int_A \sum_y
        p^{ ( d)}_{T-t}(y-x) \lambda(y,t) dx = o\left(
        m^{-t/(2+\varepsilon)} \right).
        \label{psi-p}
\end{equation}
\end{fact}

\begin{fact} {\rm (R\'ev\'esz \cite{revesz}, (6.11))}
\label{f:lambda}
        There exists a constant $C=C(m,d) >0$ such that for all
        $1\le t<T$,
\begin{equation}
       \be\left( \, \sum_{y\in \br^d} \left\{ \lambda(y,T) -
       \be\left[ \lambda(y,T)\, | \, {\cal F}(t)\right] \right\}^2 \right)
       \le C \, {m^{2T - t} \over  (T-t)^{d/2}} .
       \label{lambda-lambda}
\end{equation}
\end{fact}

\medskip

\begin{lemma}
\label{l:branching}
       Let $\varepsilon>0$. Almost surely for all large $t$, we have
       $\lambda(y,t)=0$ whenever $\|y\|>t^{1+\varepsilon}$.
\end{lemma}

\medskip

\noindent {\it Proof.} This follows from the usual estimate for the tail of
the Wiener process, the Borel--Cantelli lemma, and
(\ref{N0}).\hfill$\diamondsuit$

\medskip

\begin{lemma}
\label{l:martingale}
        Let $\alpha\in Z_+^d$, and let
\begin{equation}
       V_\alpha(t) =\sum_y H_\alpha (y, t) \lambda(y,t).
       \label{mm2.3}
\end{equation}
       Then, $({1\over m^t}\, V_\alpha(t), \; t\ge 0)$ is a martingale and
$$ N_\alpha := \lim_{t\to \infty} {V_\alpha (t) \over m^t}
$$
\noindent exists and is finite almost surely.
\end{lemma}

\medskip

\noindent {\it Proof.} We start by proving the martingale property. Recall
that $\psi (\br^d,t)$ stands for the total number of particles at time
$t$. Thus, by numbering these particles and considering them all starting
from time $t=0$ (many of them share common paths, at least partially),
we can write $\sum_y H_\alpha (y, t) \lambda(y,t) = \sum_{i=1}^{\psi
(\br^d,t)} H_\alpha (W^{ (  i)}(t), t)$, where $(W^{ ( i)}, \; i\ge 1)$ is a
sequence of $\br^d$-valued Wiener processes (they are {\it not}
independent). Conditioning on
${\cal F}(t-1)$ and on $\psi (\br^d,t)$, we have
\begin{eqnarray*}
      &&\be \left( \sum_{i=1}^{\psi (\br^d,t)} H_\alpha (W^{ (  i)}(t),  t) \,
\Big|
\, {\cal F}(t-1), \; \psi (\br^d,t) \right)
        \\
      &=& \sum_{i=1}^{\psi (\br^d,t)} H_\alpha (W^{ ( i)}(t-1), t-1)
        \\
      &=& \sum_{i=1}^{\psi (\br^d,t-1)} Y_{i,t-1} H_\alpha (W^{ (  i)}(t-1),
t-1),
\end{eqnarray*}

\noindent the last identity following from the fact that many particles at
time $t$ come from the same ancestor at time $t-1$, with $Y_{i,t-1}$
denoting the number of offspring from the $i$-th particle at time $t-1$.

Integrating on both sides gives that
\begin{eqnarray}
      &&\be \left( \sum_{i=1}^{\psi (\br^d,t)} H_\alpha (W^{ (  i)}(t),  t) \,
\Big|
\, {\cal F}(t-1) \right)\label{ui.0}
        \\
      &=& \sum_{i=1}^{\psi (\br^d,t-1)}  \be(Y) H_\alpha (W^{ ( i)}(t-1),
       t-1)\nonumber
        \\
      &=& m\sum_{i=1}^{\psi (\br^d,t-1)} H_\alpha (W^{ ( i)}(t-1), t-1),
\nonumber
\end{eqnarray}

\noindent proving that $t\mapsto {1\over m^t}  V_\alpha(t)$ is a
martingale.

We now show that $({1\over m^t}\, V_\alpha(t), \; t\ge 0)$ converges to
a finite limit almost surely. With the above notation we first write
\begin{equation} V_\alpha(t)=\sum_{l=1}^{\psi (\br^d,t)} H_\alpha
(W^{ ( l)}(t), t)=\sum_{l=1}^{\psi (\br^d,t-1)}\sum_{m=1}^{Y_{l,t-1}}
H_\alpha(W^{ ( l,m)}(t), t)\label{ui.9}
\end{equation} where $W^{ ( l,m)}(t)$ is the $m$-th child of the $l$-th
particle which dies at time $t-1$. Then we can write
\begin{eqnarray} V_\alpha(t)^{ 2}&=&\sum_{l=1}^{\psi
(\br^d,t-1)}\sum_{m=1}^{Y_{l,t-1}} H^{ 2}_\alpha(W^{ ( l,m)}(t),
t)\label{ui.10}\\ &+&
\sum_{l=1}^{\psi (\br^d,t-1)}\sum_{m\neq
n,\,m,n=1}^{Y_{l,t-1}}H_\alpha(W^{ ( l,m)}(t), t)H_\alpha(W^{ ( l,n)}(t),
t)\nonumber\\ &+&
\sum_{i\neq j,\,i,j=1}^{\psi
(\br^d,t-1)}\sum_{m=1}^{Y_{i,t-1}}\sum_{n=1}^{Y_{j,t-1}}H_\alpha(W^{
( i,m)}(t), t)H_\alpha(W^{ ( j,n)}(t), t)\nonumber
\end{eqnarray} Therefore
\begin{eqnarray} &&
\be \left( V_\alpha(t)^{ 2}
\, \Big|
\, {\cal F}(t-1), \; \psi (\br^d,t) \right)\label{ui.1}
        \\
      &=&\sum_{i=1}^{\psi (\br^d,t-1)}Y_{i,t-1}  \be \left(  H^{ 2}_\alpha
(W^{ (  i,1)}(t), t)
\, \Big|
\, {\cal F}(t-1) \right)\nonumber
        \\ &+&
\sum_{i=1}^{\psi (\br^d,t-1)}(Y^{ 2}_{i,t-1}-Y_{i,t-1}) H^{
2}_\alpha(W^{ ( i)}(t-1),t-1)\nonumber\\
      &+& \sum_{i\neq j,\,i,j=1}^{\psi (\br^d,t-1)} Y_{i,t-1}Y_{j,t-1}
    H_\alpha (W^{ (  i)}(t-1),t-1)H_\alpha (W^{ (  j)}(t-1),t-1).
    \nonumber
\end{eqnarray}

Thus
\begin{eqnarray} &&
\be \left( V_\alpha(t)^{ 2}
\, \Big|
\, {\cal F}(t-1) \right)\label{ui.2}
        \\
      &=&\sum_{i=1}^{\psi (\br^d,t-1)}m  \be \left(  H^{ 2}_\alpha (W^{ (
i)}(t), t)
\, \Big|
\, {\cal F}(t-1) \right)\nonumber
        \\&+&
\sum_{i=1}^{\psi (\br^d,t-1)}(\si^{ 2}+m^{ 2}-m) H^{ 2}_\alpha(W^{ (
i)}(t-1),t-1)\nonumber\\
      &+& \sum_{i\neq j,\,i,j=1}^{\psi (\br^d,t-1)} m^{ 2}
    H_\alpha (W^{ (  i)}(t-1),t-1)H_\alpha (W^{ (  j)}(t-1),t-1).
    \nonumber\\
    &=&  \sum_{i=1}^{\psi (\br^d,t-1)} \left[  m
       \be \left( H_\alpha^2 (W^{ ( i)}(t), t) \, \Big| \, {\cal
       F}(t-1) \right) \right. \nonumber
       \\
    && \qquad\qquad \qquad \left. + (\si^{ 2}-m)  H_\alpha^2
       (W^{ ( i)}(t-1), t-1) \right] +m^2 V_\alpha(t-1)^{ 2}.\nonumber
\end{eqnarray}

Recall that $\be ( \psi (\br^d,t-1)) = m^{t-1}$ (Athreya and Ney
\cite{athreya-ney}, p.~9). It is easy to see using (\ref{3m.2}) that
$\be\left(
       H_\alpha^2 (W^{(1)}(t), t) \right)=\al !t^{ |\al|}$. Hence
\begin{eqnarray} &&
\be \left( V_\alpha(t)^{ 2} \right)\label{ui.11}\\ && = m^{t-1}
\al ! (mt^{ |\al|}+(\si^{ 2}-m)( t-1)^{ |\al|}) +m^2\be \left(
V_\alpha(t-1)^{ 2} \right)\nonumber\\ && = m^{t-1}
\al ! (m(t^{ |\al|}-( t-1)^{ |\al|})+\si^{ 2}( t-1)^{ |\al|}) +m^2\be \left(
V_\alpha(t-1)^{ 2} \right).\nonumber
\end{eqnarray} This gives us that
\begin{equation} 0<\be \left( {V_\alpha(t)^{ 2} \over m^{
2t}}-{V_\alpha(t-1)^{ 2} \over m^{ 2( t-1)}}\right)\leq c{t^{ |\al|} \over
m^{ t}}.\label{ui.12}
\end{equation} Hence, using the fact that $V_\alpha(t)/m^{ t}$ is a
martingale we have that
\begin{eqnarray} \be \left( \sum_{ t=1}^{ \ff}\Bigg|{V_\alpha(t) \over
m^{ t}}-{V_\alpha(t-1)
\over m^{ t-1}}\Bigg| \right)&\leq &
\sum_{ t=1}^{ \ff}\left\{\be \left( \left({V_\alpha(t) \over m^{
t}}-{V_\alpha(t-1)
\over m^{ t-1}}\right)^{ 2} \right)\right\}^{ 1/2}\nonumber\\&=&
\sum_{ t=1}^{ \ff}\left\{\be \left( {V_\alpha(t)^{ 2} \over m^{ 2t}}-
{V_\alpha(t-1)^{ 2} \over m^{ 2(t-1)}} \right)\right\}^{
1/2}\label{f.3}\\&\leq & c\sum_{ t=1}^{ \ff}{t^{ |\al|/2} \over m^{ t/2}}
<\ff ,
\nonumber
\end{eqnarray} so that
\begin{equation}
\sum_{ t=1}^{ \ff}\Bigg|{V_\alpha(t) \over m^{ t}}-{V_\alpha(t-1)
\over m^{ t-1}}\Bigg|<\ff, \qquad \hbox{\rm a.s.}\label{ui.13}
\end{equation} This shows that $({1\over m^t}\, V_\alpha(t), \; t\ge 0)$
converges to a finite limit almost  surely.\hfill$\diamondsuit$
\medskip

\begin{remark} {\rm  Note that  by induction from (\ref{ui.11})
\begin{equation}\qquad
\be \left( V_\alpha(t)^{ 2} \right)= m^{t-1}\al ! (
\si^{ 2}\sum_{ j=1}^{ t-1}m^{t- j}j^{ |\al|}+m\sum_{ j=1}^{ t}m^{t- j}(j^{
|\al|}-( j-1)^{ |\al|}))\label{ui.6}
\end{equation} and therefore
\begin{equation}\qquad
\be \left( N_\alpha ^{ 2} \right)= m^{-1}\al ! (
\si^{ 2}\sum_{ j=1}^{\ff}m^{- j}j^{ |\al|}+m \sum_{ j=1}^{
\ff}m^{- j}(j^{ |\al|}-( j-1)^{ |\al|})).\label{f.2}
\end{equation} } 
\end{remark}

\begin{lemma}
\label{l:N}
       Let $\alpha\in Z_+^d$, and let $V_\alpha$, $N_\alpha$ be as in
Lemma \ref{l:martingale}. Then for any
       $\varepsilon>0$, we have that almost surely as $t\to \infty$,
\begin{equation}
       {V_\alpha (t) \over m^t} = N_\alpha + o\left( m^{-
       t/(2+\varepsilon)}\right) .
       \label{N}
\end{equation}
\end{lemma}

\medskip

\noindent {\it Proof.}  We claim that
\begin{equation} {V_\alpha (t^{ 2}) \over m^{ t^{ 2}}} = \be\left(
{V_\alpha (t^{ 2}) \over m^{ t^{ 2}}}
\,\Bigg |\, {\cal F}(t)\right) + o\left( m^{- t/(2+2\varepsilon)}\right) ,
\qquad
\hbox{\rm a.s.}\label{c.1}
\end{equation}

To see this, we first observe that by Fact \ref{f:lambda}, Chebyshev's
inequality and the Borel--Cantelli lemma that almost surely for $t\to
\infty$,
$$
\max_{y\in \br^d} \left| \lambda(y,t^2) - \be\left( \lambda(y,t^2)\, | \,
{\cal F}(t)\right) \right| = o\left( m^{t^2- t/(2+\varepsilon)}\right) .
$$

\noindent Assembling this estimate with (\ref{mm2.3}) and Lemma
\ref{l:branching}, together with the fact that
$\sup_{ \|y\|\leq t^{2(1+\varepsilon)} }H_{ \al}(y,t^2)\le c
t^{2(1+\varepsilon) |\alpha|}$, we get (\ref{c.1}).

Since $\be ( {V_\alpha (t^2) \over m^{t^2}} \, | \, {\cal F}(t))=
{V_\alpha(t) \over m^{t}}$ (by Lemma \ref{l:martingale}), it follows
from (\ref{c.1}) that
$$  {V_\alpha (t^2) \over m^{t^2}} - {V_\alpha (t) \over m^t} = o\left(
m^{- t/(2 + 2\varepsilon)}\right) ,\qquad \hbox{\rm a.s.}
$$

\noindent As a consequence,
\begin{equation}\qquad N_\alpha -{V_\alpha (t) \over m^t}=\sum_{
j=0}^{ \ff}\left( {V_\alpha (t^{ 2^{ j+1}}) \over m^{ t^{ 2^{j+1}}}}
    - {V_\alpha (t^{ 2^{ j}}) \over m^{ t^{ 2^{ j}}}}\right) = o\left( m^{-
t/(2 + 2\varepsilon)}\right) ,\qquad \hbox{\rm a.s.}\label{con.1}
\end{equation}

\noindent This proves our lemma, since $\varepsilon>0$ is
arbitrary.\hfill$\diamondsuit$

\bigskip

We have now all the ingredients to prove Theorem \ref{t:main}.

\bigskip

\noindent {\it Proof of Theorem \ref{t:main}.} Fix $k\ge 1$. Fix
$0<\gamma<{1\over 2( k+1)}$, and let $t= \lfloor T^\gamma\rfloor$.
Let
$\varepsilon>0$ be such that $(1+\varepsilon)\gamma <{1\over 2(
k+1)}$.  We will show that, almost surely for $T
\to
\infty$,
\begin{eqnarray}
        {\psi(A,T) \over m^{T}}
     &=& {1\over (2\pi T)^{d/2}} \sum_{n=0}^k {(-T)^{-n}
        \over 2^{n} }\sum_{|\alpha|=n}{ 1\over \al !}  \sum_{\bb\preceq
2\al} { 2\al\choose\bb }  ( -1)^{ |\bb|} M_{\bb}( A)   {V_{2\al-\bb}( t)
\over m^{t}}
\nn\\ &&+ o\left( T^{-(k+d/2)}\right) + O\left(
        m^{-t/(2+\varepsilon)}\right) ,\label{fin.1}
\end{eqnarray}

\noindent where $V_{2\al-\bb}$ is defined in (\ref{mm2.3}). Our
Theorem will then follow from Lemma \ref{l:N}.

By Fact \ref{f:cond-exp}, we have, almost surely for $T \to \infty$,
$$
     {\psi(A,T) \over m^{T}} =   {1 \over m^{t}} \int_A \sum_y p^{ (
d)}_{T-t}(y-x)
\lambda(y,t) dx + o\left( m^{- t/(2+\varepsilon)}\right).
$$ On the other hand we can write
\begin{equation}\qquad
    (2\pi T)^{d/2}p^{ ( d)}_{T-t}(y-x) = {1\over (1-t/T )^{d/2}} \exp\left( -
{\|y-x\|^2\over 2 (T-t )}\right)=f(z,t,x,y) ,
\label{fin.2}
\end{equation} where $z=1/T$ and
\begin{equation} f(z,t,x,y)={1\over (1-tz )^{d/2}} \exp\left( -
{\|y-x\|^2z\over 2 (1-tz )}\right)\label{fin.3}
\end{equation}
    is a $C^{ \ff}$ function of $z$ near $z=0$ as long as $tz \ll 1$. If we
expand $f(z,t,x,y)$ in a finite Taylor series in $z$ around $z=0$, it is
clear that we can bound the remainder $R_{ k+1}(z,t,x,y)$ of order
$k+1$ by a polynomial in $\|y-x\|$ of order at most $2(k+1 )$.

According to Lemma \ref{l:branching}, almost surely for all large
$T$, $\lambda(y,t)=0$ as long as $\|y\|>T^{(1+\varepsilon)\gamma}$.
Together with (\ref{N0}) which implies that the number of points $y$
with
$\la(y,t)\neq 0$ is bounded by $cm^{ t}$ and the fact that  $A$ is
bounded we have
\begin{equation} {1 \over m^{t}} \int_A \sum_y R_{ k+1}(T^{ -1},t,x,y)
\lambda(y,t) dx\leq cT^{ 2(1+\varepsilon)\ga ( k+1)}=o(T).\label{fin.4}
\end{equation}

By inspection of Lemma
\ref{ld:kernel}, the first $k$ terms in the Taylor series for $f(z,t,x,y)$
give rise to the the first line of (\ref{fin.1}), completing the proof of  that
formula and hence of our Theorem.
\hfill$\diamondsuit$

\section{$L^{ p}$ convergence}\label{s:Lp}

In this section we show that if the offspring distribution $Y$ has $p$
moments for some even integer $p$ then ${V_\alpha(t) \over m^{ t}}$
converges in $L^{ p}$.

Introduce the notation
   \[\wt{\sum}_{ i_{ 1},\ldots,i_{ j}=1}^{ n} =:\sum_{\stackrel{ i_{
1},\ldots,i_{ j}=1}{ i_{ l}\neq i_{ m},\,\forall l\neq m}}^{ n}\] for
summation over non-repeated indices. Let $Z_{ t}=\psi (\br^d,t)$,
$F_{\al; \,i }( t)=H_{\al }(W^{ (i)}(t), t)$ and
\begin{equation} U_{\al^{ ( 1)},\ldots, \al^{ ( p)}}( t)=\wt{\sum}_{ i_{
1},\ldots,i_{p}=1}^{ Z_{ t}}\prod_{ h=1}^{p} F_{\al^{ ( h)}; \,i_{ h} }( t).
\label{4.13}
\end{equation}

The following Lemma will play an important role in showing that
${V_\alpha(t) \over m^{ t}}$ converges in $L^{ p}$.

\begin{lemma}\label{lem-pbound} Let $k$ be an integer with $\be ( |Y|^{
k})<\ff$. Then for any $\al^{ ( 1)},\ldots, \al^{ ( k)}$ we can find
$c,\bb<\ff$ independent of $t$ such that
\begin{equation}
\left |\be \left (U_{\al^{ ( 1)},\ldots, \al^{ ( k)}}( t) \right)\right |
\leq ct^{ \bb }m^{ kt}.\label{4.14}
\end{equation}
\end{lemma}

\noindent {\it Proof of Lemma \ref{lem-pbound}.} We will prove this
Lemma by induction on $k$. The case of $k=1$ is trivial. Assume that we
have proven this Lemma for all $k\leq p-1$.

We can write
\begin{eqnarray} && U_{\al^{ ( 1)},\ldots, \al^{ ( p)}}( t)=\wt{\sum}_{
i_{ 1},\ldots,i_{p}=1}^{ Z_{ t}}\prod_{ h=1}^{p} F_{\al^{ ( h)}; \,i_{ h} }(
t)\label{4.50}\\ && = \sum_{ k=1}^{ p}\wt{\sum}_{ i_{
1},\ldots,i_{k}=1}^{ Z_{ t-1}}\sum_{ A_{ 1}\cup\cdots\cup A_{
k}=[1,p]}\prod_{ h=1}^{ k}
\left ( \wt{\sum}_{j_{s}=1,\forall s\in A_{ h}}^{Y_{ i_{ h},t-1}}
\prod_{ m\in A_{ h}} F_{\al^{ (m)};\,i_{ h},j_{ m} }( t)
\right)\nonumber
\end{eqnarray} where the sum $\sum_{ A_{ 1}\cup\cdots\cup A_{
k}=[1,p]}$ runs over all partitions of $[1,p]=\{ 1,\ldots,p\}$ by
$k$ non-empty sets $A_{ 1},\ldots, A_{ k}$ and
$F_{\al;\,l,m }( t)=H_\alpha(W^{ ( l,m)}(t), t)$. Introducing the falling
factorial notation $(x)_{ k}=x( x-1)\cdots ( x-k+1))$  we have that
\begin{eqnarray}&&
\be \left (\prod_{ h=1}^{ k}
\left ( \wt{\sum}_{j_{s}=1,\forall s\in A_{ h}}^{Y_{ i_{ h},t-1}}
\prod_{ m\in A_{ h}} F_{\al^{ (m)};\,i_{ h},j_{ m} }( t)
\right)\, \Big|
\, {\cal F}(t-1) \right) \label{4.51}\\ &&=
\prod_{ h=1}^{ k} \be \left (  \left ( Y\right)_{ |A_{ h}|}\right)
\prod_{ m\in A_{ h}} F_{\al^{ (m)};\,i_{ h} }( t-1).
\nonumber
\end{eqnarray} Hence
\begin{eqnarray}&&
\be \left (U_{\al^{ ( 1)},\ldots, \al^{ ( p)}}( t)
\, \Big|
\, {\cal F}(t-1) \right) \label{4.52}\\ &&=\sum_{ k=1}^{ p}\wt{\sum}_{
i_{ 1},\ldots,i_{k}=1}^{ Z_{ t-1}}\sum_{ A_{ 1}\cup\cdots\cup A_{
k}=[1,p]}
\prod_{ h=1}^{ k} \be \left (  \left ( Y\right)_{ |A_{ h}|}\right)
\prod_{ m\in A_{ h}} F_{\al^{ (m)};\,i_{ h} }( t-1)\nonumber\\ &&=m^{
p}U_{\al^{ ( 1)},\ldots, \al^{ ( p)}}( t-1) \nonumber\\ &&\hspace{
.3in}+\sum_{ k=1}^{ p-1}\wt{\sum}_{ i_{ 1},\ldots,i_{k}=1}^{ Z_{
t-1}}\sum_{ A_{ 1}\cup\cdots\cup A_{ k}=[1,p]}
\prod_{ h=1}^{ k} \be \left (  \left ( Y\right)_{ |A_{ h}|}\right)
\prod_{ m\in A_{ h}} F_{\al^{ (m)};\,i_{ h} }( t-1).
\nonumber
\end{eqnarray}

Note that by (\ref{3m.2})
\begin{eqnarray} &&
\sum_{ n=0}^{ \ff}{r^{ n} \over n!}H_{ n}( x,t)\sum_{ m=0}^{ \ff}{s^{ m}
\over m!}H_{ m}( x,t)=e^{ r x-tr^{ 2}/2}e^{ s x-ts^{ 2}/2}.\label{4.15}\\
&&=e^{ (r+s) x-t(r+s)^{ 2}/2}e^{ trs}\nonumber\\ &&=\sum_{ j=0}^{
\ff}{(r+s)^{ j} \over j!}H_{ j}( x,t)\sum_{ k=0}^{
\ff}{( trs)^{ k}\over k!}\nonumber\\ &&=\sum_{ j=0}^{ \ff}\sum_{
i=0}^{j} {r^{ i}s^{ j-i}\over i!(j-i)!}H_{ j}( x,t)\sum_{ k=0}^{
\ff}{( trs)^{ k}\over k!}.\nonumber
\end{eqnarray} Equating coefficients of $r^{ n}s^{ m}$ we find that
\begin{equation} H_{ n}( x,t)H_{ m}( x,t)=n!m!\sum_{ k=0}^{m\wedge n}
{t^{ k}\over k!}{1 \over (n-k )!(m-k )!}H_{ n+m-2k}( x,t).\label{4.16}
\end{equation} Using this to reduce products of Hermite functions to
sums we find that
\begin{eqnarray} &&
\be \left (U_{\al^{ ( 1)},\ldots, \al^{ ( p)}}( t)
\, \Big|
\, {\cal F}(t-1) \right)=m^{ p}U_{\al^{ ( 1)},\ldots, \al^{ ( p)}}(
t-1)\label{4.53}\\ && \hspace{ 1in}+
\sum_{ j=1}^{ p-1}\sum_{ \bb^{ ( 1)},\ldots, \bb^{ ( j)}} c(\al;\, p;\bb^{
( 1)},\ldots, \bb^{ ( j)};t ) U_{\bb^{ ( 1)},\ldots, \bb^{ ( j)}}(
t-1)\nonumber
\end{eqnarray} where $\sum_{ \bb^{ ( 1)},\ldots, \bb^{ ( j)}}$ is a finite
sum over
$ \bb^{ ( 1)},\ldots, \bb^{ ( j)}$ such that
   $\sum_{ l=1}^{ j} |\bb^{( l)}|\leq \sum_{ l=1}^{ p} |\al^{( l)}|$ and the
$c(\al;\, p;\bb^{ ( 1)},\ldots, \bb^{ ( j)};t )$ are polynomials in $t$.
Hence by our induction hypothesis
\begin{equation}
\be \left (U_{\al^{ ( 1)},\ldots, \al^{ ( p)}}( t)\right)= m^{ p}\be \left
(U_{\al^{ ( 1)},\ldots, \al^{ ( p)}}( t-1)\right) +\mathcal{R}_{\al^{ (
1)},\ldots, \al^{ ( p)}}( t)\label{4.54}
\end{equation} with $|\mathcal{R}_{\al^{ ( 1)},\ldots, \al^{ ( p)}}( t)|\leq
ct^{ \bb}m^{ ( p-1)( t-1)}$ for some $\bb,c<\ff$ independent of $t$.
Iterating this completes the proof of our Lemma for
$k=p$.\hfill$\diamondsuit$

   \begin{prop}\label{lem-pconv} Let $p$ be an even integer with $\be (
|Y|^{ p})<\ff$. Then ${V_\alpha(t) \over m^{ t}}$ converges in $L^{ p}$.
\end{prop}

\noindent {\it Proof of Proposition \ref{lem-pconv}.} Note that because
of the presence of the polynomial factor $t^{\bb}$ in (\ref{4.14}) we
cannot simply use Lemma \ref{lem-pbound} to show that
${V_\alpha(t) \over m^{ t}}$ is bounded uniformly in $L^{ p}$. Rather,
we will show that for some $c,\bb<\ff$ independent of $t$
\begin{equation}
\left| \be \left ( \left\{     V_\alpha(t)- mV_\alpha(t-1)\right\}^{p}
\right ) \right| \leq ct^{ \bb}m^{ t( p-1)}.\label{4.17}
\end{equation} Then
\begin{equation}
\left| \be \left ( \left\{     {V_\alpha(t)\over m^{ t}}- {V_\alpha(t-1)\over
m^{ t-1}}\right\}    ^{ p} \right ) \right|
\leq ct^{ \bb }m^{ -t}\label{4.18}
\end{equation} and therefore (it is here that we need $p$ even)
\begin{equation}
   \sum_{ t=1}^{ \ff}\left\|    {V_\alpha(t)\over m^{ t}}-
{V_\alpha(t-1)\over m^{ t-1}}\right\|_{ p}\leq c\sum_{ t=1}^{\ff} t^{
\bb /p}m^{ -t/p}<\ff\label{4.19}
\end{equation} which will complete the proof of the proposition.

The basic idea of the proof of (\ref{4.17}) is that the subtraction
eliminates the highest order term in the expectation leaving only sums of
terms of the form $U_{\al^{ ( 1)},\ldots, \al^{ ( k)}}( t)$ with $k\leq p-1$.

We now prove (\ref{4.17}). We have that
\begin{eqnarray}&&
\be \left ( \left\{     V_\alpha(t)- mV_\alpha(t-1)\right\}    ^{ p} \right
)\label{4.2}\\ &&=\sum_{ k=0}^{ p}{p\choose k}( -1)^{ k}m^{ k}
\be \left (     V^{ p-k}_\alpha(t) V^{ k}_\alpha(t-1)  \right )\nonumber\\
&&=\sum_{ k=0}^{ p}{p\choose k}( -1)^{ k}m^{ k}
\be \left (     \be \left (V^{ p-k}_\alpha(t)\, \Big|
\, {\cal F}(t-1) \right)    V^{k}_\alpha(t-1)
\right).\nonumber
\end{eqnarray}

By (\ref{ui.9}) we have
\begin{equation} V_\alpha(t) =\sum_{l=1}^{Z_{
t-1}}\sum_{m=1}^{Y_{l,t-1}} F_{\al;\,l,m }( t)\label{4.1}
\end{equation} where  $F_{\al;\,l,m }( t)=H_\alpha(W^{ ( l,m)}(t), t)$.
Thus
\begin{eqnarray}&& V^{n}_\alpha(t)\label{4.3}\\ &&=\sum_{j=1}^{ n}
\wt{\sum}_{ i_{ 1},\ldots,i_{ j}=1}^{ Z_{ t-1}}\sum_{ l_{ 1}+\cdots+l_{
j}=n}{n\choose l_{ 1},\ldots,l_{ j}}\prod_{ h=1}^{ j} \left (
\sum_{r=1}^{Y_{ i_{ h},t-1}} F_{\al;\,i_{ h},r }( t)   \right)^{ l_{ h}}
\nonumber\\ &&=
\wt{\sum}_{ i_{ 1},\ldots,i_{ n}=1}^{ Z_{ t-1}}\prod_{ h=1}^{ n} \left (
\sum_{r=1}^{Y_{ i_{ h},t-1}} F_{\al;\,i_{ h},r }( t)   \right)
\nonumber\\ &&+\sum_{j=1}^{ n-1}
\wt{\sum}_{ i_{ 1},\ldots,i_{ j}=1}^{ Z_{ t-1}}\sum_{ l_{ 1}+\cdots+l_{
j}=n}{n\choose l_{ 1},\ldots,l_{ j}}\prod_{ h=1}^{ j} \left (
\sum_{r=1}^{Y_{ i_{ h},t-1}} F_{\al;\,i_{ h},r }( t)   \right)^{ l_{ h}}
\nonumber
\end{eqnarray} and
\begin{eqnarray}&&
\left ( \sum_{r=1}^{Y_{ i_{ h},t-1}} F_{\al;\,i_{ h},r }( t) 
\right)^{ l_{ h}}
\label{4.4}\\ &&=
\sum_{s=1}^{ l_{ h}}
\wt{\sum}_{ r_{ 1},\ldots,r_{ s}=1}^{ Y_{ i_{ h},t-1}}\sum_{ q_{
1}+\cdots+q_{ s}=l_{ h}}{l_{ h}\choose q_{ 1},\ldots,q_{ s}}\prod_{
f=1}^{ s} F^{ q_{ f}}_{\al;\,i_{ h},r_{ f} }( t).
\nonumber
\end{eqnarray}

Thus
\begin{eqnarray}&&
\be \left ( \left ( \sum_{r=1}^{Y_{ i_{ h},t-1}} F_{\al;\,i_{ h},r }( t)
\right)^{ l_{ h}}\, \Big|
\, {\cal F}(t-1) \right) \label{4.5}\\ &&=
\sum_{s=1}^{ l_{ h}}\be \left (  \left ( Y\right)_{ s}\right)
\sum_{ q_{ 1}+\cdots+q_{ s}=l_{ h}}{l_{ h}\choose q_{ 1},\ldots,q_{
s}}\prod_{ f=1}^{ s}
\be \left ( F^{ q_{ f}}_{\al;\,i_{ h} }( t)\, \Big|
\, {\cal F}(t-1) \right).
\nonumber
\end{eqnarray}

Using (\ref{4.16}) to reduce products of Hermite functions to sums we
find that by (\ref{4.3})-(\ref{4.5}) we can write, with $\al^{
(i)}=\al,\,i=1,\ldots,n$
\begin{eqnarray} &&
\be \left (  V^{n}_\alpha(t)\, \Big|
\, {\cal F}(t-1) \right)=m^{ n}U_{\al^{ ( 1)},\ldots, \al^{ ( n)}}(
t-1)\label{4.30}\\  && \hspace{ 1in}+
\sum_{ j=1}^{ n-1}\sum_{ \bb^{ ( 1)},\ldots, \bb^{ ( j)}} c(\al;\, n;\bb^{
( 1)},\ldots, \bb^{ ( j)};t ) U_{\bb^{ ( 1)},\ldots, \bb^{ ( j)}}(
t-1)\nonumber
\end{eqnarray} where $\sum_{ \bb^{ ( 1)},\ldots, \bb^{ ( j)}}$ is a finite
sum and the
$c(\al;\, n;\bb^{ ( 1)},\ldots, \bb^{ ( j)};t )$ are polynomials in $t$.

We next observe that
\begin{eqnarray}&& V^{n}_\alpha(t-1)=\left (  \sum_{l=1}^{Z_{ t-1}}
F_{\al;\,l }( t)  \right)^{ n} \label{4.8}\\ &&=\sum_{j=1}^{ n}
\wt{\sum}_{ i_{ 1},\ldots,i_{ j}=1}^{ Z_{ t-1}}\sum_{ l_{ 1}+\cdots+l_{
j}=n}{n\choose l_{ 1},\ldots,l_{ j}}\prod_{ h=1}^{ j} F^{ l_{ h}}_{\al;\,i_{
h}}( t-1)
\nonumber\\ &&=U_{\al^{ ( 1)},\ldots, \al^{ ( n)}}( t-1)
\nonumber\\ &&\hspace{ 1in}+\sum_{j=1}^{ n-1}
\wt{\sum}_{ i_{ 1},\ldots,i_{ j}=1}^{ Z_{ t-1}}\sum_{ l_{ 1}+\cdots+l_{
j}=n}{n\choose l_{ 1},\ldots,l_{ j}}\prod_{ h=1}^{ j} F^{l_{ h}}_{\al;\, i_{
h}}( t-1)
\nonumber\\ &&=U_{\al^{ ( 1)},\ldots, \al^{ ( n)}}( t-1)
\nonumber\\ &&\hspace{ 1in}+\sum_{ j=1}^{ n-1}
\sum_{ \ga^{ ( 1)},\ldots, \ga^{ ( j)}} d(\al;\, n;\ga^{ ( 1)},\ldots, \ga^{ (
j)};t ) U_{\ga^{ ( 1)},\ldots, \ga^{ ( j)}}( t-1)
\nonumber
\end{eqnarray} where we have again used (\ref{4.16}) to reduce
products of Hermite functions to sums, and the $d(\al;\, n;\ga^{ (
1)},\ldots, \ga^{ ( j)};t )$ are polynomials in $t$.

Similarly
\begin{eqnarray} && U_{\bb^{ ( 1)},\ldots, \bb^{ ( j)}}( t-1)U_{\ga^{ (
1)},\ldots, \ga^{ ( k)}}( t-1)\label{4.31}\\ && = \left ( \wt{\sum}_{ i_{
1},\ldots,i_{j}=1}^{ Z_{ t-1}}\prod_{ h=1}^{j} F_{\bb^{ (h)};\,i_{ h} }(
t-1)  \right)     \left (\wt{\sum}_{j_{ 1},\ldots,j_{k}=1}^{ Z_{ t-1}}\prod_{
l=1}^{k} F_{\ga^{ ( l)};\,j_{ l} }( t-1)\right) \nonumber\\ &&=U_{\bb^{ (
1)},\ldots, \bb^{ ( j)}, \ga^{ ( 1)},\ldots, \ga^{ ( k)}}( t-1)
\nonumber\\ &&\hspace{ .5in}+\sum_{ m=1}^{ j+k-1}
\sum_{ \ze^{ ( 1)},\ldots, \ze^{ ( m)}} f(\bb, \ga;\ze^{ ( 1)},\ldots, \ze^{
( m)};t ) U_{\ze^{ ( 1)},\ldots, \ze^{ ( m)}}( t-1)
\nonumber
\end{eqnarray} where we have abreviated $\bb=(\bb^{ ( 1)},\ldots,
\bb^{ ( j)}),\,\ga=(
\ga^{ ( 1)},\ldots, \ga^{ ( k)} )$.

Combining (\ref{4.30})-(\ref{4.31}) we have that for each $k\leq p$
\begin{eqnarray} &&
    m^{ k} \be \left (V^{ p-k}_\alpha(t)\, \Big|
\, {\cal F}(t-1) \right)    V^{k}_\alpha(t-1)
\label{4.32}\\ && =m^{ p}U_{\al^{ ( 1)},\ldots, \al^{ ( p)}}(t-1)
\nonumber\\ &&\hspace{ .7in}+\sum_{ j=1}^{ p-1}
\sum_{ \ga^{ ( 1)},\ldots, \ga^{ ( j)}} f(\al;\, n;\ga^{ ( 1)},\ldots, \ga^{ (
j)};t ) U_{\ga^{ ( 1)},\ldots, \ga^{ ( j)}}( t-1)\nonumber
\end{eqnarray} where the $f(\al;\, n;\ga^{ ( 1)},\ldots, \ga^{ ( j)};t )$ are
polynomials in $t$. Substituting back into (\ref{4.2}) and using the fact
that $\sum_{ k=0}^{ p}{p\choose k}( -1)^{ k}=0$ we find that the
$m^{ p}U_{\al^{ ( 1)},\ldots, \al^{ ( p)}}(t-1)  $'s cancel, and we can write
\begin{eqnarray} &&
\be \left ( \left\{     V_\alpha(t)- mV_\alpha(t-1)\right\}^{p} \right
)\label{4.33}\\ && =\sum_{ j=1}^{ p-1}
\sum_{ \ga^{ ( 1)},\ldots, \ga^{ ( j)}} g(\al;\, n;\ga^{ ( 1)},\ldots, \ga^{ (
j)};t )
\be \left (U_{\ga^{ ( 1)},\ldots, \ga^{ ( j)}}( t-1)\right)\nonumber
\end{eqnarray} where the $g(\al;\, n;\ga^{ ( 1)},\ldots, \ga^{ ( j)};t )$ are
polynomials in $t$. (\ref{4.14}) then completes the proof of (\ref{4.17})
and hence of our Proposition.
\hfill$\diamondsuit$

\begin{remark} {\rm  Note that by Proposition \ref{lem-pconv} we have
that $\|{V_\alpha(t)\over m^{ t}}\|_{ p}$ is bounded uniformly in $t$, so
that
\begin{equation}
   \left\| V_\alpha(t)   \right\|_{ p}\leq cm^{ t}.\label{4.80}
\end{equation} Arguing as before, any $U_{\al^{ ( 1)},\ldots, \al^{ ( p)}}(
t)$, where
   $\al^{ ( 1)},\ldots, \al^{( k)}$ are now arbitrary, can be written as
\begin{equation} U_{\al^{ ( 1)},\ldots, \al^{ ( p)}}( t)=\prod_{ i=1}^{
p}V_{\al^{ ( i)} }(t)+\mbox{ terms of `lower order'}\label{4.81}
\end{equation} and thus using (\ref{4.80}), H\"older's inequality and
(\ref{4.14}) for $k\leq p-1$ we can refine (\ref{4.14}) and find
$c,\bb<\ff$ independent of $t$ such that
\begin{equation}
\left |\be \left (U_{\al^{ ( 1)},\ldots, \al^{ ( p)}}( t) \right)\right |
\leq cm^{ pt}.\label{4.14a}
\end{equation} (Here we require that $Y$ have $r$ momnets for some
even $r\geq p$). } 
\end{remark}

\bigskip

{\footnotesize

\baselineskip=12pt

\hskip20pt P\'al R\'ev\'esz

\hskip20pt Institut f\"ur Statistik und Wahrscheinlichkeitstheorie

\hskip20pt Technische Universit\"at Wien

\hskip20pt Wiedner Hauptstrasse 8-10/107

\hskip20pt A-1040 Vienna

\hskip20pt Austria

\hskip20pt {\tt revesz@ci.tuwien.ac.at}

\bigskip
\bigskip

\noindent
\begin{tabular}{lll} & \hskip20pt Jay Rosen
     & \hskip20pt Zhan Shi \\  & \hskip20pt Department of Mathematics
     & \hskip20pt Laboratoire de Probabilit\'es UMR 7599 \\    &
\hskip20pt College of Staten Island, CUNY
     & \hskip20pt Universit\'e Paris VI \\    & \hskip20pt Staten Island, NY
10314
     & \hskip20pt 4 place Jussieu \\    & \hskip20pt U.S.A.
     & \hskip20pt F-75252 Paris Cedex 05 \\   & \hskip20pt {\tt
jrosen3@earthlink.net}
     & \hskip20pt France \\  &
     & \hskip20pt zhan@proba.jussieu.fr
\end{tabular}

}

\end{document}